\documentclass[11pt]{article}
\usepackage[russian]{babel}
\usepackage[cp1251]{inputenc}
\usepackage{amsmath,latexsym}
\usepackage[psamsfonts]{amssymb}
\usepackage[dvips]{graphicx}
\usepackage[mathcal]{euscript}
\begin{document}
MSC2010.35C06
\begin{center}
\textbf{BUILDING SINGULAR SOLUTIONS
FOR DEGENERATE HIGH ORDER EQUATIONS\\}
\textbf{B.Y.Irgashev}

Namangan Engineering-Construction Institute, Namangan, 160103, Uzbekistan\\
e-mail: bahromirgasev@gmail.com\\
http://orcid.org/0000-0001-7204-9127
\end{center}

\textbf{Abstract:}
In this paper, using the similarity method, we construct particular solutions with singularities for degenerate high-order equations. The considered equations have singularities of the first and second kind. Particular solutions are expressed through generalized hypergeometric functions. Sufficient conditions are obtained under which the number of linearly independent solutions is equal to the order of the equation.


In our works [1], [2], particular solutions were constructed with singularities for high-order equations with multiple characteristics with constant coefficients. In this paper, we will construct self-similar solutions for equations of the form
\[{x^\alpha }D_x^pu\left( {x,y} \right) - {y^\beta }D_y^qu\left( {x,y} \right) = 0,\,\,\,p > q,\,\,0 \le \alpha  < p,\,0 \le \beta  < q,\eqno(1)\]
\[{y^\beta }D_x^pu\left( {x,y} \right) - {x^\alpha }D_y^qu\left( {x,y} \right) = 0,\,p > q,\,\,\alpha ,\beta  \ge 0,\eqno(2)\]
\[{x^\alpha }{y^\beta }D_x^pu\left( {x,y} \right) - D_y^qu\left( {x,y} \right) = 0,\,\,\,p > q,\,\,0 \le \alpha  < p,\beta  \ge 0,\eqno(3)\]
\[D_x^pu\left( {x,y} \right) - {x^\alpha }{y^\beta }D_y^qu\left( {x,y} \right) = 0,\,\,\,p > q,\,\,\alpha  \ge 0,\,0 \le \beta  < q.\eqno(4)\]
Note that in (3), a similarity method was used to study a third-order equation that has first-order degeneracy.

We will search for self-similar solutions of equations (1) - (4) in the form
\[u\left( {x,t} \right) = {y^b}v\left( t \right),\,\,t = x{y^a},\]
parameters $a, b$ are to be determined. We have
\[D_x^pu = {y^{b + ap}}D_t^pv(t),\eqno(5)\]
\[
D_y^q u\left( {x,y} \right) = \sum\limits_{k = 0}^q {C_q^k \left(
{y^b } \right)^{\left( k \right)} } D_y^{q - k} v\left( t
\right),\eqno(6)
\]
\[
\begin{array}{l}
 D_y^j v\left( t \right) = D_y^{j - 1} \left( {D_t^1 vD_y^1 t} \right) = \sum\limits_{k_1  = 0}^{j - 1} {\left( {C_{j - 1}^{k_1 } D_y^{j - 1 - k_1 } \left( {D_y^1 t} \right)D_y^{k_1 } \left( {D_t^1 v} \right)} \right)}  =  \\
  = D_y^{j - 1} \left( {D_y^1 t} \right)D_t^1 v + \sum\limits_{k_1  = 1}^{j - 1} {C_{j - 1}^{k_1 } D_y^{j - 1 - k_1 } \left( {D_y^1 t} \right)D_y^{k_1  - 1} \left( {D_t^2 vD_y^1 t} \right)} . \\
 \end{array}
\]

We introduce the following notation $\overline{\left(a\right)}_{j+
1}=a\left({a - 1}\right)\left({a - 2}\right)...\left({a-
j}\right)$ and assume that $\overline{\left(a\right)}_0=
1,$ then $D_y^j\left({D_y^1t}\right)=\frac{{\overline
{\left(a\right)}_{j+1}t}}{{y^{j+1}}}$. Hence we have
\[
\begin{array}{l}
 D_y^j v\left( t \right) = \frac{{\overline {\left( a \right)} _j t}}{{y^j }}D_t^1 v(t) + \sum\limits_{k_1  = 1}^{j - 1} {C_{j - 1}^{k_1 } \frac{{\overline {\left( a \right)} _{j - k_1 } }}{{y^{j - k_1 } }}t} \sum\limits_{k_2  = 0}^{k_1  - 1} {C_{k_1  - 1}^{k_2 } } D_y^{k_1  - k_2  - 1} \left( {D_y^1 t} \right)D_y^{k_2 } \left( {D_t^2 v} \right) =  \\
  = \frac{{\overline {\left( a \right)} _j }}{{y^j }}tD_t^1 v + \sum\limits_{k_1  = 1}^{j - 1} {C_{j - 1}^{k_1 } \frac{{\overline {\left( a \right)} _{j - k_1 } }}{{y^{j - k_1 } }}tD_y^{k_1  - 1} \left( {D_y^1 t} \right)D_t^2 v}  +  \\
  + \sum\limits_{k_1  = 1}^{j - 1} {\sum\limits_{k_2  = 1}^{k_1  - 1} {C_{j - 1}^{k_1 } C_{k_1  - 1}^{k_2 } \frac{{\overline {\left( a \right)} _{j - k_1 } }}{{y^{j - k_1 } }}\frac{{\overline {\left( a \right)} _{k_1  - k_2 } }}{{y^{k_1  - k_2 } }}t^2 } } D_y^{k_2  - 1} \left( {D_t^3 vD_y^1 t} \right). \\
 \end{array}
 \]
Continuing this process, we obtain the following formula
\[
D_y^j v\left( {x,y} \right) = y^{ - j} \sum\limits_{i = 1}^j
{\left( {A_{i - 1}^j \left( a \right)t^i D_t^i v}
\right)},\eqno(7)
\]
where
\[
A_{i - 1}^j \left( a \right) = \sum\limits_{k_1  = i - 1}^{j - 1}
{\sum\limits_{k_2  = i - 2}^{k_1  - 1} {...} \sum\limits_{k_{i -
1}  = 1}^{k_{i - 2}  - 1} {\left( {\overline {\left( a \right)}
_{k_{i - 1} } \prod\limits_{s = 1}^{i - 1} {C_{k_{s - 1}  -
1}^{k_s } \overline {\left( a \right)} _{k_{s - 1}  - k_s } } }
\right)} } ,
\]
moreover, $k_0  = j,\,\,i - 1 = \overline {1,j - 1} ,\,\,k_1  > k_2
> ... > k_{i - 1}  \ge 1.$\\
We note some properties of the coefficients $A_i^{j}$

\textbf{\emph{Lemma.}}\\
  \emph{
1.$ A_i^{i + 1} \left( a \right) = a^{i + 1} ;$\\
2.$ A_i^j \left( a \right) = a\left( {\left( {i + 1} \right)A_i^{j
- 1} \left( a \right) + A_{i - 1}^{j - 1} \left( a \right)}
\right) - \left( {j - 1} \right)A_i^{j - 1} \left( a \right); $\\
3.$ A_0^j \left( a \right) = \overline {\left( a \right)_j };$ \\
4.$ \sum\limits_{j = 1}^s {\overline {\left( x \right)_j } A_{j -
1}^s \left( a \right)}  = \overline {\left( {ax} \right)_s };$}\\
\emph{5. $\sum\limits_{s = 0}^n {C_n^s } \overline {\left( b \right)}
_{n - s} \overline {\left( y \right)} _s  = \overline {\left( {y +
b} \right)} _n.$}\\

\emph{\textbf{Proof.}} \\
 1). Obviously, $ A_0 ^ 1 = a $
 using property 2 of this
lemma and the fact that $ A_i ^ j = 0 $ for $ i \ge j$
,we have
 \[
A_i^{i + 1}  = a\left( {\left( {i + 1} \right)A_i^i  + A_{i - 1}^i
} \right) - iA_i^i  = aA_{i - 1}^i  = a^i A_0^1  = a^{i + 1} .
\]
2). We have
\[
D_y^{j - 1} u = y^{ - j + 1} \sum\limits_{i = 1}^{j - 1} {\left(
{A_{i - 1}^{j - 1} t^i D_t^i u} \right)} ,
\]
we take the derivative
\[
\left( {D_y^{j - 1} u} \right)_y  = D_y^j u = \left( {1 - j}
\right)y^{ - j} \sum\limits_{i = 1}^{j - 1} {\left( {A_{i - 1}^{j
- 1} t^i D_t^i u} \right)}  + y^{ - j + 1} \sum\limits_{i = 1}^{j
- 1} {\left( {iA_{i - 1}^{j - 1} t^{i - 1} \frac{{at}}{y}D_t^i u}
\right)}  +
\]
\[
 + y^{ - j + 1} \sum\limits_{i = 1}^{j - 1} {\left( {A_{i - 1}^{j - 1} t^i \frac{{at}}{y}D_t^{i + 1} u} \right)}  = y^{ - j} \sum\limits_{i = 1}^{j - 1} {t^i D_t^i u\left( {\left( {1 - j} \right)A_{i - 1}^{j - 1}  + iA_{i - 1}^{j - 1} a} \right)}  +
\]
\[
 + y^{ - j} \sum\limits_{i = 2}^j {\left( {aA_{i - 2}^{j - 1} t^i D_t^i u} \right)} ,
\]
on the other side
\[
D_y^j u = y^{ - j} \sum\limits_{i = 1}^j {\left( {A_{i - 1}^j t^i
D_t^i u} \right)} ,
\]
equating the coefficients at the same degrees, we obtain
required result. \\
3).  from relation (2) automatically follows
for $ j = 0. $ \\
4). We use property 2 of the coefficients $ A_i ^ j $ and assuming that
$ A_ {- 1} ^ {s - 1} = 0 $, we have
\[
\sum\limits_{j = 1}^s {\overline {\left( x \right)_j } A_{j - 1}^s
}  = \sum\limits_{j = 1}^s {\overline {\left( x \right)_j } }
\left( {a\left( {jA_{j - 1}^{s - 1}  + A_{j - 2}^{s - 1} } \right)
- \left( {s - 1} \right)A_{j - 1}^{s - 1} } \right) =
\]
$$
 = a\sum\limits_{j = 1}^s {\overline {\left( x \right)_j } \left( {jA_{j - 1}^{s - 1}  + A_{j - 2}^{s - 1} } \right)}  - \left( {s - 1} \right)\sum\limits_{j = 1}^{s - 1} {\overline {\left( x \right)_j } A_{j - 1}^{s - 1} }= $$
 $$= a\sum\limits_{j = 1}^{s - 1} {A_{j - 1}^{s - 1} } \left( {j\overline {\left( x \right)_j }  + \overline {\left( x \right)_{j + 1} } } \right) - \left( {s - 1} \right)\sum\limits_{j = 1}^{s - 1} {\overline {\left( x \right)_j } } A_{j - 1}^{s - 1}  =
$$
\[
 = ax\sum\limits_{j = 1}^{s - 1} {\overline {\left( x \right)_j } A_{j - 1}^{s - 1} }  - \left( {s - 1} \right)\sum\limits_{j = 1}^{s - 1} {\overline {\left( x \right)_j } A_{j - 1}^{s - 1} }  = \left( {ax + 1 - s} \right)\sum\limits_{j = 1}^{s - 1} {\overline {\left( x \right)_j } } A_{j - 1}^{s - 1} ,
\]
Continuing the process we obtain
\[
\sum\limits_{j = 1}^s {\overline {\left( x \right)_j } A_{j - 1}^s
\left( a \right)}  = \overline {\left( {ax} \right)_s } .
\]
5). Shown in the work of M. Aigner [4] \\
\emph{\textbf{Lemma is proved.}} \\
Putting (7) in (6) we have
\[
D_y^q u\left( {x,y} \right) = \sum\limits_{k = 0}^{q - 1} {C_q^k
\overline {\left( b \right)_k } y^{b - k} } y^{ - q + k}
\sum\limits_{j = 1}^{q - k} {\left( {A_{j - 1}^{q - k} \left( a
\right)t^j D_t^j v\left( t \right)} \right)}  + \overline {\left(
b \right)_q } y^{b - q} v\left( t \right) =
\]
\[
 = y^{b - q} \left( {\sum\limits_{k = 0}^{q - 1} {C_q^k \overline {\left( b \right)_k } } \sum\limits_{j = 1}^{q - k} {\left( {A_{j - 1}^{q - k} \left( a \right)t^j D_t^j v\left( t \right)} \right)}  + \overline {\left( b \right)_q } v\left( t \right)}
 \right).\eqno(8)
 \]
We substitute (5) and (8) into equations (1) - (4), then we will have
\[{x^\alpha }{y^{b + ap}}D_t^pv(t) - {y^{b - q + \beta }}\left( {\sum\limits_{k = 0}^{q - 1} {C_q^k\overline {{{\left( b \right)}_k}} } \sum\limits_{j = 1}^{q - k} {\left( {A_{j - 1}^{q - k}\left( a \right){t^j}D_t^jv\left( t \right)} \right)}  + \overline {{{\left( b \right)}_q}} v\left( t \right)} \right) = 0,\]
\[{y^{b + ap + \beta }}D_t^pv(t) - {x^\alpha }{y^{b - q}}\left( {\sum\limits_{k = 0}^{q - 1} {C_q^k\overline {{{\left( b \right)}_k}} } \sum\limits_{j = 1}^{q - k} {\left( {A_{j - 1}^{q - k}\left( a \right){t^j}D_t^jv\left( t \right)} \right)}  + \overline {{{\left( b \right)}_q}} v\left( t \right)} \right) = 0,\]
\[{x^\alpha }{y^{b + ap + \beta }}D_t^pv(t) - {y^{b - q}}\left( {\sum\limits_{k = 0}^{q - 1} {C_q^k\overline {{{\left( b \right)}_k}} } \sum\limits_{j = 1}^{q - k} {\left( {A_{j - 1}^{q - k}\left( a \right){t^j}D_t^jv\left( t \right)} \right)}  + \overline {{{\left( b \right)}_q}} v\left( t \right)} \right) = 0,\]
\[{y^{b + ap}}D_t^pv(t) - {x^\alpha }{y^{b - q + \beta }}\left( {\sum\limits_{k = 0}^{q - 1} {C_q^k\overline {{{\left( b \right)}_k}} } \sum\limits_{j = 1}^{q - k} {\left( {A_{j - 1}^{q - k}\left( a \right){t^j}D_t^jv\left( t \right)} \right)}  + \overline {{{\left( b \right)}_q}} v\left( t \right)} \right) = 0.\]
Having made the transformations, we obtain
\[{\left( {x{y^{\frac{{ap + q - \beta }}{\alpha }}}} \right)^\alpha }D_t^pv(t) - \left( {\sum\limits_{k = 0}^{q - 1} {C_q^k\overline {{{\left( b \right)}_k}} } \sum\limits_{j = 1}^{q - k} {\left( {A_{j - 1}^{q - k}\left( a \right){t^j}D_t^jv\left( t \right)} \right)}  + \overline {{{\left( b \right)}_q}} v\left( t \right)} \right) = 0,\]
\[D_t^pv(t) - {\left( {x{y^{ - \frac{{q + ap + \beta }}{\alpha }}}} \right)^\alpha }\left( {\sum\limits_{k = 0}^{q - 1} {C_q^k\overline {{{\left( b \right)}_k}} } \sum\limits_{j = 1}^{q - k} {\left( {A_{j - 1}^{q - k}\left( a \right){t^j}D_t^jv\left( t \right)} \right)}  + \overline {{{\left( b \right)}_q}} v\left( t \right)} \right) = 0,\]
\[{\left( {x{y^{\frac{{ap + \beta  + q}}{\alpha }}}} \right)^\alpha }D_t^pv(t) - \left( {\sum\limits_{k = 0}^{q - 1} {C_q^k\overline {{{\left( b \right)}_k}} } \sum\limits_{j = 1}^{q - k} {\left( {A_{j - 1}^{q - k}\left( a \right){t^j}D_t^jv\left( t \right)} \right)}  + \overline {{{\left( b \right)}_q}} v\left( t \right)} \right) = 0,\]
\[D_t^pv(t) - {\left( {x{y^{ - \frac{{q + ap - \beta }}{\alpha }}}} \right)^\alpha }\left( {\sum\limits_{k = 0}^{q - 1} {C_q^k\overline {{{\left( b \right)}_k}} } \sum\limits_{j = 1}^{q - k} {\left( {A_{j - 1}^{q - k}\left( a \right){t^j}D_t^jv\left( t \right)} \right)}  + \overline {{{\left( b \right)}_q}} v\left( t \right)} \right) = 0.\]
Now find the parameter $a$. For the first equation we have:
\[a =  - \frac{{q - \beta }}{{p - \alpha }},\]
for the second:
\[a =  - \frac{{q + \beta }}{{p + \alpha }},\]
for the third:
\[a =  - \frac{{q + \beta }}{{p - \alpha }},\]
and finally for the fourth:
\[a =  - \frac{{q - \beta }}{{p + \alpha }}.\]
Now for the first and third equations we get
$$
{t^\alpha }D_t^pv(t) - \left( {\sum\limits_{k = 0}^{q - 1} {C_q^k\overline {{{\left( b \right)}_k}} } \sum\limits_{j = 1}^{q - k} {\left( {A_{j - 1}^{q - k}\left( a \right){t^j}D_t^jv\left( t \right)} \right)}  + \overline {{{\left( b \right)}_q}} v\left( t \right)} \right) = 0,\eqno(9)
$$
and for the second and fourth
$$
D_t^pv(t) - {t^\alpha }\left( {\sum\limits_{k = 0}^{q - 1} {C_q^k\overline {{{\left( b \right)}_k}} } \sum\limits_{j = 1}^{q - k} {\left( {A_{j - 1}^{q - k}\left( a \right){t^j}D_t^jv\left( t \right)} \right)}  + \overline {{{\left( b \right)}_q}} v\left( t \right)} \right) = 0,\eqno(10)
$$
make another replacement $z = {t^c}$, where $c$ is to be determined, we have
\[D_t^pv(t) = {t^{ - p}}\sum\limits_{j = 1}^p {A_{j - 1}^p\left( c \right){z^j}D_z^jv\left( z \right)} ,\]
\[{t^j}D_t^jv(t) = {t^j}{t^{ - j}}\sum\limits_{s = 1}^j {A_{s - 1}^j\left( c \right){z^s}D_z^sv\left( z \right)}  = \sum\limits_{s = 1}^j {A_{s - 1}^j\left( c \right){z^s}D_z^sv\left( z \right)}.\]
Substituting into equations (9) and (10), we obtain
\[\sum\limits_{j = 1}^p {A_{j - 1}^p\left( c \right){z^j}D_z^jv\left( z \right)}  - {t^{p - \alpha }}\left( {\sum\limits_{k = 0}^{q - 1} {C_q^k\overline {{{\left( b \right)}_k}} } \sum\limits_{j = 1}^{q - k} {\left( {A_{j - 1}^{q - k}\left( a \right)\sum\limits_{s = 1}^j {A_{s - 1}^j\left( c \right){z^s}D_z^sv\left( z \right)} } \right)}  + \overline {{{\left( b \right)}_q}} v\left( z \right)} \right) = 0,\]
\[\sum\limits_{j = 1}^p {A_{j - 1}^p\left( c \right){z^j}D_z^jv\left( z \right)}  - {t^{p + \alpha }}\left( {\sum\limits_{k = 0}^{q - 1} {C_q^k\overline {{{\left( b \right)}_k}} } \sum\limits_{j = 1}^{q - k} {\left( {A_{j - 1}^{q - k}\left( a \right)\sum\limits_{s = 1}^j {A_{s - 1}^j\left( c \right){z^s}D_z^sv\left( z \right)} } \right)}  + \overline {{{\left( b \right)}_q}} v\left( z \right)} \right) = 0,\]
therefore, for equation (9), the replacement has the form \[z = {t^{p - \alpha }},\] and for (10):\[z = {t^{p + \alpha }}.\]
Given this, for the study we obtain the following equation:
\[\sum\limits_{j = 1}^p {A_{j - 1}^p\left( c \right){z^j}D_z^jv\left( z \right)}  = z\left( {\sum\limits_{k = 0}^{q - 1} {C_q^k\overline {{{\left( b \right)}_k}} } \sum\limits_{j = 1}^{q - k} {\left( {A_{j - 1}^{q - k}\left( a \right)\sum\limits_{s = 1}^j {A_{s - 1}^j\left( c \right){z^s}D_z^sv\left( z \right)} } \right)}  + \overline {{{\left( b \right)}_q}} v\left( z \right)} \right).\]
For a compact record of the resulting equation, we make one more change
\[z = {e^\tau },\]
as
\[{z^j}D_z^ju\left( z \right) = {\overline {\left( D \right)} _j}u\left( \tau  \right),\]
then
\[\sum\limits_{j = 1}^p {A_{i - 1}^p\left( c \right){{\overline {\left( D \right)} }_j}v\left( \tau  \right)}  - {e^\tau }\sum\limits_{k = 0}^{q - 1} {C_q^k{{\overline {\left( b \right)} }_k}} \sum\limits_{j = 1}^{q - k} {A_{j - 1}^{q - k}\left( a \right)} \sum\limits_{s = 1}^j {\left( {A_{s - 1}^j\left( c \right){{\overline {\left( D \right)} }_s}v\left( \tau  \right)} \right)}  - {e^\tau }{\overline {\left( b \right)} _q}v\left( \tau  \right) = 0,\]
We use the fourth and fifth properties of the coefficients $A_i^j$
\[{\overline {\left( {cD} \right)} _p}v - {e^\tau }\sum\limits_{k = 0}^{q - 1} {C_q^k{{\overline {\left( b \right)} }_k}} \sum\limits_{j = 1}^{q - k} {A_{j - 1}^{q - k}\left( a \right){{\overline {\left( {cD} \right)} }_j}v}  - {e^\tau }{\overline {\left( b \right)} _q}v\left( \tau  \right) = 0,\]
from here
\[{\overline {\left( {cD} \right)} _p}v\left( \tau  \right) - {e^\tau }\sum\limits_{k = 0}^q {C_q^k{{\overline {\left( b \right)} }_k}} {\overline {\left( {acD} \right)} _{q - k}}v = 0,\]
or finally
\[{\overline {\left( {cD} \right)} _p}v\left( \tau  \right) = {e^t}{\overline {\left( {acD + b} \right)} _q}v\left( \tau  \right).\eqno(11)\]
The solution to equation (11) will be sought in the form of a series
\[v\left( \tau  \right) = \sum\limits_{n = 0}^\infty  {{c_n}{e^{\left( {n + \gamma } \right)\tau }}},\eqno(12)\]
then we have
\[{\overline {\left( {cD} \right)} _p}v\left( \tau  \right) = \sum\limits_{n = 0}^\infty  {{c_n}{{\overline {\left( {cD} \right)} }_p}{e^{\left( {n + \gamma } \right)\tau }}}  = {e^{\gamma \tau }}\sum\limits_{n = 0}^\infty  {{c_n}{{\overline {\left( {c\left( {n + \gamma } \right)} \right)} }_q}{e^{n\tau }}} ,\]
\[{e^\tau }{\overline {\left( {acD + b} \right)} _q}v\left( \tau  \right) = {e^\tau }{e^{\gamma \tau }}\sum\limits_{n = 0}^\infty  {{c_n}{{\overline {\left( {ac\left( {n + \gamma } \right) + b} \right)} }_q}{e^{n\tau }}} ,\]
from here we get
\[\sum\limits_{n = 0}^\infty  {{c_n}{{\overline {\left( {c\left( {n + \gamma } \right)} \right)} }_p}{e^{n\tau }}}  = {e^\tau }\sum\limits_{n = 0}^\infty  {{c_n}{{\overline {\left( {ac\left( {n + \gamma } \right) + b} \right)} }_q}{e^{n\tau }}} ,\]
so that on the left side there is no free coefficient, we require the equality
\[{\overline {\left( {c\gamma } \right)} _p} = 0\, \Rightarrow {\gamma _i} = \frac{i}{c}\,,\,i = 0,1,...,\left( {p - 1} \right),\]
when this condition is satisfied, we obtain
\[\sum\limits_{n = 0}^\infty  {{c_n}{{\overline {\left( {i + cn} \right)} }_p}{e^{n\tau }}}  = {e^\tau }\sum\limits_{n = 0}^\infty  {{c_n}{{\overline {\left( {ia + acn + b} \right)} }_q}{e^{n\tau }}} ,\]
or
\[\begin{array}{l}
\sum\limits_{n = 0}^\infty  {{c_n}{c^p}\left( {\frac{i}{c} + n} \right)\left( {\frac{{i - 1}}{c} + n} \right)...\left( {\frac{{i - p + 1}}{c} + n} \right){e^{n\tau }}}  = \\
 = {e^\tau }\sum\limits_{n = 0}^\infty  {{c_n}{{\left( {ac} \right)}^q}\left( {\frac{i}{c} + \frac{b}{{ac}} + n} \right)\left( {\frac{i}{c} + \frac{{b - 1}}{{ac}} + n} \right)...\left( {\frac{i}{c} + \frac{{b - q + 1}}{{ac}} + n} \right){e^{n\tau }}} ,
\end{array}\]
from here to find the coefficients $ {c_n} $ we get the formula
\[{c_n} = \frac{{{{\left( {ac} \right)}^q}\left( {\frac{i}{c} + \frac{b}{{ac}} + n - 1} \right)\left( {\frac{i}{c} + \frac{{b - 1}}{{ac}} + n - 1} \right)...\left( {\frac{i}{c} + \frac{{b - q + 1}}{{ac}} + n - 1} \right)}}{{{c^p}\left( {\frac{i}{c} + n} \right)\left( {\frac{{i - 1}}{c} + n} \right)...\left( {\frac{{i - p + 1}}{c} + n} \right)}}{c_{n - 1}} = \]
\[ = {\left( {\frac{{{a^q}}}{{{c^{p - q}}}}} \right)^n}\frac{{{{\left( {\frac{i}{c} + \frac{b}{{ac}}} \right)}_n}{{\left( {\frac{i}{c} + \frac{{b - 1}}{{ac}}} \right)}_n}...{{\left( {\frac{i}{c} + \frac{{b - q + 1}}{{ac}}} \right)}_n}}}{{{{\left( {\frac{i}{c} + 1} \right)}_n}{{\left( {\frac{{i - 1}}{c} + 1} \right)}_n}...{{\left( {\frac{{i - p + 1}}{c} + 1} \right)}_n}}}{c_0}.\]
For definiteness, let $c_0 = 1$.Then, we get the following solutions:
\[{v_i} = {e^{{\gamma _i}\tau }}\sum\limits_{n = 0}^\infty  {\frac{{{{\left( {\frac{i}{c} + \frac{b}{{ac}}} \right)}_n}{{\left( {\frac{i}{c} + \frac{{b - 1}}{{ac}}} \right)}_n}...{{\left( {\frac{i}{c} + \frac{{b - q + 1}}{{ac}}} \right)}_n}}}{{{{\left( {\frac{i}{c} + 1} \right)}_n}{{\left( {\frac{{i - 1}}{c} + 1} \right)}_n}...{{\left( {\frac{{i - p + 1}}{c} + 1} \right)}_n}}}{{\left( {\frac{{{a^q}}}{{{c^{p - q}}}}{e^\tau }} \right)}^n}}  = \]
\[ = {z^{{\gamma _i}}}\sum\limits_{n = 0}^\infty  {\frac{{{{\left( {\frac{i}{c} + \frac{b}{{ac}}} \right)}_n}{{\left( {\frac{i}{c} + \frac{{b - 1}}{{ac}}} \right)}_n}...{{\left( {\frac{i}{c} + \frac{{b - q + 1}}{{ac}}} \right)}_n}}}{{{{\left( {\frac{i}{c} + 1} \right)}_n}{{\left( {\frac{{i - 1}}{c} + 1} \right)}_n}...{{\left( {\frac{{i - p + 1}}{c} + 1} \right)}_n}}}{{\left( {\frac{{{a^q}}}{{{c^{p - q}}}}z} \right)}^n}}  = \]
\[ = {t^i}\sum\limits_{n = 0}^\infty  {\frac{{{{\left( {\frac{i}{c} + \frac{b}{{ac}}} \right)}_n}{{\left( {\frac{i}{c} + \frac{{b - 1}}{{ac}}} \right)}_n}...{{\left( {\frac{i}{c} + \frac{{b - q + 1}}{{ac}}} \right)}_n}}}{{{{\left( {\frac{i}{c} + 1} \right)}_n}{{\left( {\frac{{i - 1}}{c} + 1} \right)}_n}...{{\left( {\frac{{i - p + 1}}{c} + 1} \right)}_n}}}{{\left( {\frac{{{a^q}}}{{{c^{p - q}}}}{t^c}} \right)}^n}} ,\]
\[\,i = 0,1,...,\left( {p - 1} \right).\]
Now we write out self-similar solutions of equations (1) - (4) in terms of generalized hypergeometric functions.

For the first equation we have:
\[\begin{array}{l}
{u_i}\left( {x,y} \right) = {y^b}{\left( {x{y^{ - \frac{{q - \beta }}{{p - \alpha }}}}} \right)^i} \cdot \\
 \cdot {}_q{F_{p - 1}}\left[ {\begin{array}{*{20}{c}}
{\frac{i}{{p - \alpha }} - \frac{b}{{q - \beta }},...,\frac{i}{{p - \alpha }} - \frac{{b - q + 1}}{{q - \beta }},{{\left( { - 1} \right)}^q}\frac{{{{\left( {q - \beta } \right)}^q}}}{{{{\left( {p - \alpha } \right)}^p}}}\frac{{{x^{p - \alpha }}}}{{{y^{q - \beta }}}}}\\
{\frac{i}{{p - \alpha }} + 1,..,\frac{{i - \left( {i - 1} \right)}}{{p - \alpha }} + 1,\frac{{i - \left( {i + 1} \right)}}{{p - \alpha }} + 1,...,\frac{{i - p + 1}}{{p - \alpha }} + 1}
\end{array}} \right],
\end{array}\eqno(13)\]
\[\,i = 0,1,...,\left( {p - 1} \right),\]
here
\[
{}_pF_q \left[ {\begin{array}{*{20}c}
   {a_1 ,...,a_p ,x}  \\
   {b_1 ,...,b_q }  \\
\end{array}} \right] = \sum\limits_{k = 0}^\infty  {\frac{{\left( {a_1 } \right)_k ...\left( {a_p } \right)_k }}{{\left( {b_1 } \right)_k ...\left( {b_q } \right)_k }}} \frac{{x^k }}{{k!}}
\]
- a generalized hypergeometric function, 
\[{\left( a \right)_n} = a\left( {a + 1} \right)...\left( {a + n - 1} \right)\]
- the symbol of Pochhammer. \\
Now we choose the parameter $ b $, so that the order of the hypergeometric function decreases by one, for this we solve the following equation for $ b $:
\[\frac{i}{{p - \alpha }} - \frac{{b - q + 1}}{{q - \beta }} = \frac{{i - p + 1}}{{p - \alpha }} + 1 \Rightarrow  - \frac{{b - q + 1}}{{q - \beta }} = \frac{{ - p + 1}}{{p - \alpha }} + 1 \Rightarrow \]
\[b = q - 1 + \frac{{\left( {\alpha  - 1} \right)\left( {q - \beta } \right)}}{{p - \alpha }}.\]
We study the formula (13). For linearly independent solutions to have $ p $ pieces, it is enough to fulfill the conditions
$$
\begin{array}{l}
\alpha  \notin N,\\
\frac{i}{{p - \alpha }} - \frac{{b - s}}{{q - \beta }} \ne 0,\,\,\,i = 0,...,\left( {p - 1} \right),\,s = 0,...,\left( {q - 1} \right).
\end{array}
$$
For equation (2):
\[{u_i}\left( {x,y} \right) = {y^b}{\left( {x{y^{ - \frac{{q + \beta }}{{p + \alpha }}}}} \right)^i}{}_q{F_{p - 1}}\left[ {\begin{array}{*{20}{c}}
{\frac{i}{{p + \alpha }} - \frac{b}{{q + \beta }},...,\frac{i}{{p + \alpha }} - \frac{{b - q + 1}}{{q + \beta }},\frac{{{{\left( { - 1} \right)}^q}{{\left( {q + \beta } \right)}^q}}}{{{{\left( {p + \alpha } \right)}^p}}}\left( {\frac{{{x^{p + \alpha }}}}{{{y^{q + \beta }}}}} \right)}\\
{\frac{i}{{p + \alpha }} + 1,...,\frac{{i - \left( {i - 1} \right)}}{{p + \alpha }} + 1,\frac{{i - \left( {i + 1} \right)}}{{p + \alpha }} + 1,...,\frac{{i - p + 1}}{{p + \alpha }} + 1}
\end{array}} \right],\]
\[\,i = 0,1,...,\left( {p - 1} \right),\]
\[b = q - 1 - \frac{{\left( {\alpha  + 1} \right)\left( {q - \beta } \right)}}{{p - \alpha }}.\]
For linearly independent solutions to have $ p $ pieces, it is enough to fulfill the conditions
$$
\frac{i}{{p + \alpha }} - \frac{{b - s}}{{q + \beta }} \ne 0,\,i = 0,...,\left( {p - 1} \right),\,\,s = 0,1,...,\left( {q - 1} \right).
$$
For equation (3):
\[{u_i}\left( {x,y} \right) = {y^b}{\left( {x{y^{ - \frac{{q + \beta }}{{p - \alpha }}}}} \right)^i}{}_q{F_{p - 1}}\left[ {\begin{array}{*{20}{c}}
{\frac{i}{{p - \alpha }} - \frac{b}{{q + \beta }},...,\frac{i}{{p - \alpha }} - \frac{{b - q + 1}}{{q + \beta }},{{\left( { - 1} \right)}^q}\frac{{{{\left( {q + \beta } \right)}^q}}}{{{{\left( {p - \alpha } \right)}^p}}}\frac{{{x^{p - \alpha }}}}{{{y^{q + \beta }}}}}\\
{\frac{i}{{p - \alpha }} + 1,..,\frac{{i - \left( {i - 1} \right)}}{{p - \alpha }} + 1,\frac{{i - \left( {i + 1} \right)}}{{p - \alpha }} + 1,...,\frac{{i - p + 1}}{{p - \alpha }} + 1}
\end{array}} \right],\]
\[\,i = 0,1,...,\left( {p - 1} \right),\]
\[b = q - 1 + \frac{{\left( {\alpha  - 1} \right)\left( {q + \beta } \right)}}{{p - \alpha }}.\]
For linearly independent solutions to have $ p $ pieces, it is enough to fulfill the conditions
\[\begin{array}{l}
\,\alpha  \notin N,\\
\frac{i}{{p - \alpha }} - \frac{{b - s}}{{q + \beta }} \ne 0,\,\,\,i = 0,...,\left( {p - 1} \right),\,s = 0,...,\left( {q - 1} \right).
\end{array}\]
And finally for equation (4):
\[{u_i}\left( {x,y} \right) = {y^b}{\left( {x{y^{ - \frac{{q - \beta }}{{p + \alpha }}}}} \right)^i}{}_q{F_{p - 1}}\left[ {\begin{array}{*{20}{c}}
{\frac{i}{{p + \alpha }} - \frac{b}{{q - \beta }},...,\frac{i}{{p + \alpha }} - \frac{{b - q + 1}}{{q - \beta }},\frac{{{{\left( { - 1} \right)}^q}{{\left( {q - \beta } \right)}^q}}}{{{{\left( {p + \alpha } \right)}^p}}}\left( {\frac{{{x^{p + \alpha }}}}{{{y^{q - \beta }}}}} \right)}\\
{\frac{i}{{p + \alpha }} + 1,...,\frac{{i - \left( {i - 1} \right)}}{{p + \alpha }} + 1,\frac{{i - \left( {i + 1} \right)}}{{p + \alpha }} + 1,...,\frac{{i - p + 1}}{{p + \alpha }} + 1}
\end{array}} \right],\]
\[\,i = 0,1,...,\left( {p - 1} \right),\]
\[b = q - 1 - \frac{{\left( {\alpha  + 1} \right)\left( {q - \beta } \right)}}{{p + \alpha }}.\]
For linearly independent solutions to have $ p $ pieces, it is enough to fulfill the conditions
\[\frac{i}{{p + \alpha }} - \frac{{b - s}}{{q - \beta }} \ne 0,\,i = 0,...,\left( {p - 1} \right),\,\,s = 0,1,...,\left( {q - 1} \right).\]

\begin{center}
\textbf{References}
\end{center}

1.Irgashev, B.Y. On Partial Solutions of One Equation with Multiple Characteristics and Some Properties of the Fundamental Solution. Ukr Math J 68, 868–893 (2016). https://doi.org/10.1007/s11253-016-1263-9.

2.Irgashev, B.Y. Construction of Singular Particular Solutions Expressed via Hypergeometric Functions for an Equation with Multiple Characteristics. Diff Equat 56, 315–323 (2020). https://doi.org/10.1134/S0012266120030040

3.Hasanov, A., Ruzhansky, M. Hypergeometric Expansions of Solutions of the Degenerating Model Parabolic Equations of the Third Order. Lobachevskii J Math 41, 27–31 (2020). https://doi.org/10.1134/S1995080220010059

4.Aigner M. Combinatorial Theory. Springer-Verlag Berlin Heidelberg.1997.

\end{document}